\documentclass[12pt]{amsart}
\usepackage{epsfig}
\usepackage{amscd}
\newcommand{\ov}{\overline}
\newcommand{\un}{\underline}

\newcommand{\mat}{\ov{\cA}_t}
\newcommand{\omat}{\ov{\cA}_t \otimes \ov{\cA}_t}
\newcommand{\cmat}{\ov{\cA_t \otimes \cA}_t}
\newcommand{\qg}{{_qSL_2}}

\newcommand{\cA}{{\mathcal A}}
\newcommand{\bC}{{\hspace{-0.3pt}\mathbb C}\hspace{0.4pt}}

\newtheorem{proposition}{Proposition}

\newtheorem{lemma}{Lemma}

\title{A matrix model for Quantum $SL_2$}
\author{Charles Frohman}

\address{Department of Mathematics, University of Iowa, Iowa City, IA
52242}

\email{frohman@math.uiowa.edu}

\author{Joanna Kania-Bartoszy\'{n}ska}

\address{Department of Mathematics, Boise State University, Boise, ID
83725}

\email{kania@math.idbsu.edu}

\thanks{This research was partially supported by by NSF-DMS-9803233 and
 NSF-DMS-9971905.}
  
\begin{document}

\begin{abstract} 
We describe a topological ribbon Hopf algebra whose elements are sequences
of matrices. The algebra is a quantum version of $U(sl_2)$.
\end{abstract}

\maketitle

For each nonzero $t \in \mathbb{C}$ that is not a root of unity,
we give a quantum analog $\overline{\cA}_t$ of $U(sl_2)$. The underlying
algebra of the  model is $\prod_{n=1}^{\infty} M_n(\mathbb{C})$. Consequently, 
the algebra structure, which comes from matrix
multiplication, is independent of the variable $t$.

Define  $\cA_t$ to be the unital Hopf algebra on $X$, $Y$, $K$, $K^{-1}$,
with relations: 
\begin{align}\label{relations} KX=t^2XK, \quad KY=t^{-2}YK,\\ \label{rel}
 XY-YX= \frac{K^2-K^{-2}}{t^2-t^{-2}}, \quad KK^{-1}=1. \end{align}
The comultiplication is the algebra morphism  given by: 
\[\Delta(X)=X \otimes K + K^{-1} \otimes X,\quad  
\Delta(Y)=Y \otimes K + K^{-1} \otimes Y,\]
\[\Delta(K)=K \otimes K.\]
The antipode is the antimorphism given by $S(X)=-t^2X$, $S(Y)=-t^{-2}Y$, $S(K)=K^{-1}$, and the counit
is the morphism given by $\epsilon(X)=\epsilon(Y)=0$, and $\epsilon(K)=1$.

The standard representations $\un{m}$, where $m$ is a nonnegative integer,
of $\cA_t$ have basis $e_i$, where $i$ runs
in integer steps from $-m/2$ to $m/2$. Hence as a vector space $\un{m}$ has dimension $m+1$. Recall that
\[ [n]=\frac{t^{2n}-t^{-2n}}{t^2-t^{-2}},\]
and  $[n]!=[n][n-1]\ldots[1]$.

The action of $\cA_t$ is given by 
\begin{align*}
X \cdot e_i & = [m/2 + i + 1]e_{i+1} \quad \text{but} \quad X \cdot
e_{m/2} = 0, \\ 
Y \cdot e_i & = [m/2 - i + 1]e_{i-1} \quad \text{but}
\quad Y \cdot e_{-m/2} = 0, \\
K \cdot e_i & = t^{2i}e_i.
\end{align*}

The representation $\un{m}$ can be seen as a homomorphism 
\[\rho_m: A_t \rightarrow
M_{m+1}(\bC).\] 

\begin{lemma}\label{onto} The homomorphisms $\rho_m: A_t \rightarrow
M_{m+1}(\bC)$ are onto.\end{lemma}

\proof Using the ordered basis, $\{e_{-m/2},\ldots e_{m/2}\}$,
$\rho_m(X)$ is the matrix that is zero except on the first subdiagonal,
where the entries going from the top to the bottom are $1,[2],[3],\ldots,
[m]$. Similarly, the matrix $\rho_m(Y)$ is zero except on the first 
superdiagonal, where starting from the bottom and going up the entries
are $1,[2],[3],\ldots,[m]$. 
The image of  $X^nY^p$  is a matrix with zero entries except on a particular
super- or sub-diagonal , whose distance from the diagonal is
$|n-p|$. Starting from the top,
the first $\text{min}\{p,n\}$ entries  of that diagonal 
are zero, and the subsequent entries are all nonzero. 
Thus there exist
linear combinations of the matrices $\rho_m(X^nY^p)$, with $p\leq n$,
corresponding to each of the elementary matrices whose only nonzero entry
lies on the $n-p$ subdiagonal, or on the diagonal. We are using the
pattern of zero and nonzero entries on the $n-p$ subdiagonal to see this.
By a similar analysis of $\rho_m(Y^pX^n)$ we see that 
all elementary matrices where the nonzero entry
lies on a superdiagonal can be written as a linear combination of
the  $\rho_m(Y^pX^n)$.
Since $M_{m+1}(\mathbb{C})$ is spanned by the elementary matrices, this
finishes the proof. 
\qed

Define the linear functionals ${^mc}^i_j : A_t \rightarrow \bC$  by
letting ${^mc}^i_j(Z)$ be the $ij$-th coefficient of the matrix $\rho_m(Z)$.
Let $\qg$ be the stable subalgebra of the Hopf algebra dual $A_t^o$ 
generated by linear functionals
${^mc}^i_j$. 

\begin{proposition}
The linear functionals ${^mc}^i_j$ form a basis for the algebra $\qg$.
\end{proposition}

\proof
Since \[\un{m}\otimes \un{n}=\bigoplus_{q=|m-n|}^{m+n}\un{q},\]
the linear functionals ${^mc}^i_j$ span the algebra $\qg$.
We need to show that they are also linearly independent.
The quantum Casimir is given by
\begin{equation}\label{casimir}
C =\frac{
( tK - t^{-1}K^{-1})^2
}
{(t^2-t^{-2})^2} 
+ YX \in \cA_t.
\end{equation}
 Since $C$ is central in $\cA_t$, it acts as scalar 
multiplication in
any irreducible representation. In fact, it acts on $\un{m}$ as 
$\lambda_m=\frac{(t^{m+1}-t^{-m-1})^2}{(t^2-t^{-2})^2}$. Let 
\begin{equation}\label{cmn}
C_{m,n}=\frac{C-\lambda_n}{\lambda_m-\lambda_n}.
\end{equation}
 Notice that $C_{m,n}$ 
is zero
under $\rho_n$ and is sent to the identity in $\rho_m$. The product
\begin{equation}\label{dmN}
 D_{m,N}=\prod_{p=1, p\neq m}^{N} C_{m,n}
\end{equation} 
is an element of $\cA_t$ that is 
sent to $0$ in all of the representations from $\un{1}$ to $\un{N}$, 
except $\un{m}$ where it is sent to the identity matrix.

If
some linear combination $\sum \alpha_{i,j,n}{^nc}^i_j$ is equal to zero,
it means that for all $Z \in \cA_t$,
\[ \sum_{i,j,n} \alpha_{i,j,n}{^nc}^i_j(Z)=0.\]
Let $N$ be the largest $n$ for which $\alpha_{i,j,n}\neq 0$.
For each $m$ such that $\alpha_{i,j,m}\neq 0$, 
apply the functional $\sum_{i,j,n} \alpha_{i,j,n}{^nc}^i_j$ to 
$D_{m,N}Z$. Since
\[ \sum_{i,j,n} \alpha_{i,j,n}{^nc}^i_j(D_{m,N}Z)=0,\]
it follows that for all $Z \in \cA_t$, and fixed $m$,
\[  \sum_{i,j} \alpha_{i,j,m}{^mc}^i_j(Z)=0.\]
Finally, from  lemma \ref{onto}  the homomorphisms $\rho_m$ are surjective, so the 
independence of the  ${^mc}^i_j$
follows from the independence of the matrix coefficients on $M_{m+1}(\bC)$.
Therefore all the $\alpha_{i,j,m}=0$.
\qed

The product of any two matrix coefficients can be written as a linear
combination of matrix  coefficients
\begin{equation}\label{gammas}
 {^mc}^i_j\cdot {^nc}^k_l= \sum_{u,v,p}
 \gamma_{u,v,p}^{i,j,m,k,l,n}(t)\ {^pc}^u_v
\end{equation}

Since the functionals ${^pc}^u_v$ are linearly independent, the coefficients
 $\gamma_{u,v,p}^{i,j,m,k,l,n}(t)$ are uniquely defined. 
The $\gamma_{u,v,p}^{i,j,m,k,l,n}(t)$ are versions of the Clebsch-Gordan
coefficients.
Notice that  $|m-n| \leq p \leq m+n$, 
consequently for each tuple
$(i,j,m,k,l,n)$ there are only finitely many $(u,v,p)$ with 
$\gamma_{u,v,p}^{i,j,m,k,l,n}(t) \neq 0$.

A similar computation can be performed with the analogously defined
${^mc}^i_j$ associated to $Sl_2(\mathbb{C})$. The limit as $t$ approaches 
$1$ of the coefficients
$\gamma_{u,v,p}^{i,j,m,k,l,n}(t)$ gives the corresponding quantities for
$Sl_2(\mathbb{C})$.

Let \[ \mat = M_1(\bC) \times M_2(\bC) \times M_3(\bC) \times \ldots\]
be the Cartesian product of all the matrix algebras over $\bC$ given the
product topology.  
\begin{proposition}
The homomorphism 
\begin{equation}\label{theta}
\Theta : \cA_t \rightarrow \mat
\end{equation}
given by $\Theta(Z)=(\rho_0(Z),\rho_1(Z),\rho_2(Z), \ldots)$ is
injective and its image is dense in $\mat$.
\end{proposition}

\proof
The fact that the $\rho_m$ are onto and the existence of the elements
$C_{m,n}$ defined by equation (\ref{cmn}) can
be used to prove that the image of $\Theta$ is dense in $\mat$. 

A version
of the Poincar\'{e}-Birkhoff-Witt theorem says that the monomials
$K^mX^nY^p$ form a basis for $\cA_t$ as a vector space. Using the
relation $ XY-YX= \frac{K^2-K^{-2}}{t^2-t^{-2}}$, this can be replaced
by the basis $Z_{m,n,p}$, with $Z_{m,n,p}=K^mX^nY^p$ for $n\geq p$ and
$Z_{m,n,p}=K^mY^pX^n$, when $n<p$.
In order to prove that the map $\Theta$ is injective,
consider an element $\sum\alpha_i Z_{m_i,n_i,p_i}\in\cA_t$. It
is our goal to show that if $\Theta(\sum\alpha_i Z_{m_i,n_i,p_i})=0$
then all $\alpha_i$ are zero.

In any representation the image of $Z_{m,n,p}$ is a matrix
that is zero off of the super (or sub)-diagonal corresponding to
$n-p$. Thus it suffices to consider the  sums where $n_i-p_i$ is a constant,
as long as we only work with the parts of the matrices in the image
that lie on the super- or sub- diagonal corresponding to that constant.

Assume that $n_i\geq p_i$, 
The argument is similar when $n_i <p_i$.
Suppose that, for $k\geq 0$, the image under $\Theta$ of
\begin{equation}\label{sum}\sum\alpha_i K^{m_i}X^{p_i+k}Y^{p_i},\end{equation}
on the $k$th subdiagonal is zero.
The map $\Theta$  takes $ K^{m}X^{p+k}Y^{p}$ 
to a sequence of matrices such that the 
first $p$ entries along the $k$-th subdiagonal are zero.  
Let $p$ be the minimum of the $p_i$ appearing in (\ref{sum}).
The  $(p +1)$-st entry of each $k$-th subdiagonal of each matrix
in the sequence $\Theta(\sum\alpha_i K^{m_i}X^{p_i+k}Y^{p_i})$ is  
 the image under $\Theta$ of the  collection of terms in (\ref{sum})
with $p_i=p$. All the other terms are mapped to matrices with a  
zero there.
Thus it is enough to show that whenever all the $(p+1)$-st entries on the $k$-th subdiagonal
in each entry 
of  $\Theta(\sum\alpha_i K^{m_i}X^{p+k}Y^p)$  
are zero, then  all $\alpha_i$ are zero.

Assume that all the $(p+1)$-st entries on the $k$-th subdiagonal
of $\Theta(\sum\alpha_i Z_{m_i,p+k,p})$ are zero. 
Make a sequence consisting of the $(p+1)$-st entries of the $k$-th diagonal
of the image of $Z_{0,p+k,p}$. This sequence is:
\[(0,0,\ldots,[p+k]!\prod_{r=1}^p[k+r],[p+k]!\prod_{r=1}^p[k+r+1],\ldots),\]
where the first nonzero entry corresponds to the representation
$\rho_{p+k+1}$. 
 Hence, the sequence corresponding to $Z_{m_i,p+k,p}$ is
\[(0,0,\ldots,t^{m_i(p+k)}[p+k]!\prod_{r=1}^p[k+r],
t^{m_i(p+k-1)}[p+k]!\prod_{r=1}^p[k+r+1],\ldots).\]
Supposing that we have $J$ terms in our sum, we can truncate
these sequences to get a $J\times J$ matrix, so that the coefficients
$\alpha_i$ as a column vector, must be in the kernel of that matrix.
Notice that the coefficient of the power of $t$ in each column is
the same product of quantized integers. Hence its
determinant is a product of quantized integers times the
determinant of the matrix,
\[ \begin{pmatrix} t^{m_1(p+k)} & t^{m_1(p+k-1)}& \ldots
\\ t^{m_2(p+k)} & t^{m_2(p+k-1)}& \ldots \\
\vdots & \vdots &\vdots \end{pmatrix}.\]

Factoring out a large power of $t$ from each row we get the Vandermonde
determinant,
\[\left| \begin{matrix} 1 & t^{-m_1} & t^{-2m_1} & \ldots
\\ 1 & t^{-m_2}& t^{-2m_2} &\ldots \\
\vdots & \vdots &\vdots & \vdots  \end{matrix}\right| ,\]
 which is nonzero
as long as the $t^{m_i}$ are not equal to one another. Since $t$ was
chosen specifically not to be a root of unity, all the $\alpha_i$ must
be zero. \qed

The topology induced on $\mat$ by its image under $\Theta$ is the weak
topology from $\qg$. That is a sequence $Z_n$ is Cauchy 
if for every $\phi \in \qg$, $\phi(Z_n)$ is a Cauchy sequence of 
complex
numbers. Hence $\mat$ is the completion of $\cA_t$ by equivalence classes of
Cauchy sequences in the weak topology from $\qg$.

Let $e_{i,j}(m) \in \mat$ be the sequence of matrices that is the zero
matrix in every entry, except the $m+1$-st, where it is the elementary matrix 
that is all
zeroes except for a $1$ in the $ij$-th entry. Notice that the 
$e_{i,j}(m)$ are dual
to the $^mc^i_j$ 
in the sense that $^mc^i_j(e_{k,l}(p))$ is zero 
unless the
indices are identical, in which case it is one. Also notice that any 
$A \in \mat$
can be written uniquely as $\sum_{i,j,m}\alpha_{i,j,m}e_{i,j}(m)$. The 
infinite sum makes sense!
\begin{proposition}
The algebra $\mat$ has a structure of a topological ribbon Hopf algebra.
\end{proposition}
\proof
We need to define comultiplication on $\mat$.
Every element of $\omat$ can be written as an infinite sum,
\begin{equation}\label{tensor}
 \sum_{i,j,m,k,l,n}\tau_{i,j,m,k,l,n} e_{i,j}(m)\otimes e_{k,l}(n)
\end{equation}
so that no $e_{i,j}(m)\otimes e_{k,l}(n)$ is repeated.
There are infinite sums of this form that cannot be 
decomposed as a finite sum of 
tensors of elements
of $\mat$. We topologize 
$\omat$ by saying that a sequence $W_n$ is Cauchy
if and only if for every ${^mc^i_j}\otimes {^nc^k_l}$ the sequence 
$\left({^mc^i_j}\otimes {^nc^k_l}\right)(W_n)$ is Cauchy. 
Let $\cmat$ be the completion of $\omat$ 
by
equivalence classes of Cauchy sequences. Notice that every sum of 
the
type like in equation (\ref{tensor})
yields an equivalence class of Cauchy sequences in $\omat$ by truncating to
get a sequence of partial sums.  Conversely,
if $Z_n \in \omat$ is Cauchy, by applying the 
${^mc^i_j}\otimes {^nc^k_l}$ to the sequence, and taking the limit
we get the coefficients of a unique expression of the type
(\ref{tensor}),
and two Cauchy sequences
are equivalent if and only if they give rise to the same expression.
Hence we can identify $\cmat$ with the set of expressions like in
equation (\ref{tensor}).

In order to define the comultiplication on $\mat$ with  values
in $\cmat$,   take the
adjoint of multiplication on $\qg$. Use $<\ , \ >$ to denote evaluation
of elements of $\mat$ on $\qg$, and  extend this to evaluating elements
of $\qg \otimes \qg$ on elements of $\cmat$ pairwise. Then,
\[<{^mc^i_j}\otimes {^nc^k_l},\Delta(e_{u,v}(q))>
= <{^mc^i_j}\cdot{^nc^k_l},e_{u,v}(q)>
=\gamma_{i,j,m,k,l,n}^{u,v,q}.\]

Therefore,

\[ \Delta(e_{u,v}(q))
=\sum_{i,j,m,k,l,n} \gamma_{i,j,m,k,l,n}^{u,v,q} e_{i,j}(m)\otimes 
e_{k,l}(n).\]
The sum makes sense for an arbitrary element of $\mat$ as there are 
only finitely
many nonzero $\gamma_{i,j,m,k,l,n}^{u,v,q}$ for any 
$e_{i,j}(m)\otimes e_{k,l}(n)$.
So one can sum 
\[ \Delta(\sum_{i,j,m}\alpha_{u,v,p}e_{u,v}(q))=\sum \alpha_{i,j,m}\Delta(e_{i,j}(m))=\]
\[ \sum \alpha_{i,j,m}\gamma_{i,j,m,k,l,n}^{u,v,q}(t) e_{i,j}(m)\otimes e_{k,l}(n).\]
Comultiplication is continuous  since its composition with every
${^mc^i_j}\otimes {^nc^k_l}$ is continuous.

Let $q=t^4$.
The standard formula for the universal
$R$-matrix \cite{Ka} in the Jimbo-Drinfeld model of $U_h(sl_2)$ is
\begin{equation}\label{R}
 R=\sum_{n\geq 0}\frac{(q-q^{-1})^n}{[n]}q^{-n(n+1)/2}t^{H\otimes H +nH\otimes 1 - 1\otimes nH}(X^n \otimes Y^n). 
\end{equation}

Recall that the standard Drinfeld-Jimbo model \cite{Ka} of
$U_h(sl_2)$ is generated by  $X,Y,H$. If we let $K=t^H$ then the
relations (\ref{relations}), (\ref{rel}) for $\cA_t$ can be derived from the
relations for the Drinfeld-Jimbo model. 
Consequently, interpret $H$ as the traditional image of $H$
under the standard irreducible representations of $U(sl_2)$. That is,
$H$ is the sequence of matrices, 
\[(1,\begin{pmatrix} -1 & 0\\ 0 & 1 \end{pmatrix},\begin{pmatrix} -2 & 0 & 0\\
 0 & 0 & 0 \\ 0 & 0 & 2 \end{pmatrix},...).\] 
Taking $t$ raised to this sequence gives the sequence
 $\Theta(K)$, 
where $\Theta$ is defined in equation (\ref{theta}). 
Interpret $X$ and $Y$ as the sequences of matrices coming
from the standard representations of $\cA_t$, i.e., $\Theta(X)$ and
$\Theta(Y)$. 
The resulting
expression (\ref{R}) makes sense as an element of $\omat$ since
in any particular irreducible representation only finitely many
terms are nonzero. Thus the $R$ matrix is well defined as an element of
$\omat$, and has the desired properties. \qed

\end{document}